\documentclass[11pt]{article}
\usepackage{graphicx}

\def\be{\begin{equation}}
\def\ee{\end{equation}}
\def\ff#1{\mbox{\boldmath $#1$} }

\def\e{\epsilon}
\def\x{\ff{x}}

\begin{document}

\title{Firefly Algorithm: Recent Advances and Applications}

\author{Xin-She Yang \\
School of Science and Technology, \\
Middlesex University, The Burroughs, London NW4 4BT, UK. \\
\and
Xingshi He \\
School of Science, Xi'an Polytechnic University, \\
No. 19 Jinhua South Road, Xi'an 710048, P. R. China. }

\date{}

\maketitle
\begin{abstract}
Nature-inspired metaheuristic algorithms, especially those based on swarm intelligence,
have attracted much attention in the last ten years. Firefly algorithm appeared in about five years ago,
its literature has expanded dramatically with diverse applications. In this paper, we will
briefly review the fundamentals of firefly algorithm together with a selection of
recent publications. Then, we discuss the optimality associated with balancing
exploration and exploitation, which is essential for all metaheuristic algorithms.
By comparing with intermittent search strategy, we conclude that metaheuristics
such as firefly algorithm are better than the optimal
intermittent search strategy. We also analyse algorithms and their implications for
higher-dimensional optimization problems.
\end{abstract}


{\bf Reference} to this paper should be made as follows:
Xin-She Yang and Xingshi He, (2013). `Firefly Algorithm: Recent Advances and Applications',
{\it Int. J. Swarm Intelligence}, Vol. 1, No. 1, pp. 36--50. 
DOI: 10.1504/IJSI.2013.055801

\section{Introduction}

Metaheuristic algorithms form an important part of contemporary global
optimisation algorithms, computational intelligence and soft computing.
These algorithms are usually nature-inspired with multiple interacting agents.
A subset of metaheuristcs are often referred to as swarm intelligence (SI) based
algorithms, and these SI-based algorithms have been developed by mimicking
the so-called swarm intelligence charateristics of biological agents such as
birds, fish, humans and others. For example, particle swarm optimisation
was based on the swarming behaviour of birds and fish \cite{Kennedy},
while the firefly algorithm was based on the flashing pattern of tropical fireflies \cite{Yang,YangFA}
and cuckoo search algorithm was inspired by the brood parasitism of some cuckoo species \cite{YangDeb}.

In the last two decades, more than a dozen new algorithms such as particle swarm optimisation,
differential evolution, bat algorithm, firefly algorithm and cuckoo search have appeared and
they have shown great potential in solving
tough engineering optimisation problems \cite{Yang,Bansal,Floudas,Parp,YangBook,YangBA,Gandomi}.
Among these new algorithms, it has been shown that firefly algorithm is very efficient in
dealing with multimodal, global optimisation problems.

In this paper, we will first outline the fundamentals of firefly algorithm (FA), and then review the
latest developments concerning FA and its variants. We also highlight the reasons why FA is so efficient.
Furthermore, as the balance of exploration and exploitation is important to
all metaheuristic algorithms, we will then discuss the optimality related to search
landscape and algorithms. Using the intermittent search strategy and numerical experiments,
we show that firefly algorithm is significantly more efficient than intermittent
search strategy.

\section{Firefly Algorithm and Complexity}

\subsection{Firefly Algorithm}

Firefly Algorithm (FA) was first developed by Xin-She Yang in late 2007 and 2008 at Cambridge University \cite{Yang,YangFA},
which was based on the flashing patterns and behaviour of fireflies.
In essence, FA uses the following three idealized rules:

\begin{itemize}
\item Fireflies are unisex so that one firefly will be attracted to other fireflies regardless of their sex.

\item The attractiveness is proportional to the brightness, and they both decrease as their distance increases.
Thus for any two flashing fireflies, the less brighter one will move towards the brighter one.
If there is no brighter one than a particular firefly, it will move randomly.

\item The brightness of a firefly is determined by the landscape of the objective function.

\end{itemize}

As a firefly's attractiveness is proportional to the light intensity seen by adjacent fireflies,
we can now define the variation of attractiveness $\beta$ with the distance $r$ by
    \be \beta = \beta_0 e^{-\gamma r^2}, \ee
where $\beta_0$ is the attractiveness at $r=0$.

The movement of a firefly $i$ is attracted to another more attractive (brighter) firefly $j$ is determined by
\be    \x_i^{t+1} =\x_i^t + \beta_0 e^{-\gamma r^2_{ij} } (\x_j^t-\x_i^t) + \alpha_t \; \ff{\e}_i^t, \ee
where the second term is due to the attraction. The third term is randomization with $\alpha_t$ being the randomization parameter,
and $\ff{\e}_i^t$ is a vector of random numbers drawn from a Gaussian distribution or uniform distribution at time $t$.
If $\beta_0=0$, it becomes a simple random walk. On the other hand, if $\gamma=0$, it reduces to a variant of
particle swarm optimisation \cite{Yang}.
Furthermore, the randomization $\ff{\e}_i^t$ can
easily be extended to other distributions such as L\'evy flights \cite{Yang}.
A demo version of firefly algorithm implementation by Xin-She Yang, without L\'evy flights for simplicity, can
be found at Mathworks file exchange web site.\footnote{ http://www.mathworks.com/matlabcentral/fileexchange/29693-firefly-algorithm}

\subsection{Parameter Settings}

As $\alpha_t$ essentially control the randomness (or, to some extent, the diversity of solutions), we can tune this parameter during iterations so that it can vary with
 the iteration counter $t$. So a good way to
express $\alpha_t$ is to use
\be \alpha_t=\alpha_0 \delta^t, \quad 0<\delta<1), \ee
where $\alpha_0$ is the initial randomness scaling factor, and $\delta$ is essentially a cooling factor. For most applications, we can use $\delta=0.95$ to $0.97$ \cite{Yang}.

Regarding the initial $\alpha_0$, simulations show that FA will be more efficient if $\alpha_0$ is associated with the scalings of design variables.
Let $L$ be the average scale of the problem of interest, we can set
$\alpha_0=0.01L$ initially. The factor $0.01$ comes from the fact that random walks
requires a number of steps to reach the target while balancing the local exploitation
without jumping too far in a few steps \cite{YangFA,YangBook}.

The parameter $\beta$ controls the attractiveness, and parametric studies suggest that $\beta_0=1$ can be used for most applications. However, $\gamma$ should be also related to the scaling $L$. In general, we can set $\gamma=1/\sqrt{L}$. If the scaling variations are not significant, then we can set $\gamma=O(1)$.

For most applications, we can use the population size $n=15$ to $100$, though the best range is $n=25$ to $40$ \cite{Yang,YangFA}.

\subsection{Algorithm Complexity}

Almost all metaheuristic algorithms are simple in terms of complexity, and thus they are easy to implement.
Firefly algorithm has two inner loops when going through the population $n$, and one outer loop for iteration $t$.
So the complexity at the extreme case is $O(n^2 t)$. As $n$ is small (typically, $n=40$), and $t$ is large (say, $t=5000$),
the computation cost is relatively inexpensive because the algorithm complexity is linear in terms of $t$.
The main computational cost will be in the evaluations of
objective functions, especially for external black-box type objectives. This latter case is
also true for all metaheuristic algorithms. After all, for all optimisation problems, the most
computationally extensive part is objective evaluations.

If $n$ is relatively large, it is possible to use one inner loop by ranking the attractiveness
or brightness of all fireflies using sorting algorithms. In this case, the algorithm complexity
of firefly algorithm will be $O(n t \log(n))$.

\subsection{Firefly Algorithm in Applications}

Firefly algorithm has attracted much attention and has been applied to many applications \cite{Apo,Chatt,Hass,Sayadi,YangFA2010,Horng,Horng2}.
Horng et al. demonstrated that firefly-based algorithm used least computation time for digital image
compression \cite{Horng,Horng2}, while Banati and Bajaj used firefly algorithm for feature selection and showed that
firefly algorithm produced consistent and better performance in terms of time and optimality than other algorithms \cite{Banati}.

In the engineering design problems, Gandomi et al. \cite{Gandomi} and Azad and Azad \cite{Azad} confirmed that
firefly algorithm can efficiently solve highly nonlinear, multimodal design problems. Basu and Mahanti \cite{Basu}
as well as Chatterjee et al. have applied FA in antenna design optimisation and showed that FA can outperform artificial bee colony (ABC) algorithm \cite{Chatt}. In addition, Zaman and Matin have also found that FA can outperform PSO and obtained global best results \cite{Zaman}.

Sayadi et al. developed a discrete version of FA which can efficiently solve NP-hard scheduling problems \cite{Sayadi},
while a detailed analysis has demonstrated the efficiency of FA over a wide range of test problems,
including multobjective load dispatch problems \cite{Apo,YangFA,Yang12}. Furthermore, FA can also solve scheduling and travelling
salesman problem in a promising way \cite{Palit,Jati,Yousif}.

Classifications and clustering are another important area of applications of FA with excellent performance \cite{Senthil,Rajini}.
For example, Senthilnath el al. provided an extensive performance study by compared FA with
11 different algorithms and concluded that firefly algorithm can be efficiently used for clustering \cite{Senthil}.
In most cases, firefly algorithm outperform all other 11 algorithms.
In addition, firefly algorithm has also been applied to train neural networks \cite{Nandy}.

For optimisation in dynamic environments, FA can also be very efficient as shown by Farahani et al. \cite{Fara,Fara2} and Abshouri et al. \cite{Absh}.

\subsection{Variants of Firefly Algorithm}

For discrete problems and combinatorial optimisation, discrete versions of firefly algorithm have
been developed with superior performance \cite{Sayadi,Hass,Jati,Fister,Durkota}, which can be used for
travelling-salesman problems, graph colouring and other applications.

In addition, extension of firefly algorithm to multiobjective optimisation has also been investigated \cite{Apo,YangMOFA}.

A few studies show that chaos can enhance the performance of firefly algorithm \cite{Coelho,YangFAC},
while other studies have attempted to hybridize FA with other algorithms to enhance their performance \cite{Giann,Horng2,Horng3,Ramp}.

\section{Why Firefly Algorithm is So Efficient?}

As the literature about firefly algorithm expands and new variants have emerged, all pointed out the firefly algorithm can
outperform many other algorithms. Now we may ask naturally ``Why is it so efficient?".
To answer this question, let us briefly analyze the firefly algorithm itself.

FA is swarm-intelligence-based, so it has the similar advantages that other swarm-intelligence-based algorithms have.
In fact, a simple analysis of parameters suggests that some PSO variants such as Accelerated PSO \cite{YangAPSO} are
a special case of firefly algorithm when $\gamma=0$ \cite{Yang}.

However, FA has two major advantages over other algorithms: automatical subdivision and the ability of dealing with
multimodality.  First, FA is based on attraction and attractiveness decreases with distance. This leads to the fact that
the whole population can automatically subdivide into subgroups, and each group can swarm around each mode or local optimum.
Among all these modes, the best global solution can be found. Second, this subdivision allows the fireflies to be able to
find all optima simultaneously if the population size is sufficiently higher than the number of modes.
Mathematically, $1/\sqrt{\gamma}$ controls the average distance of a group of fireflies that can be seen by
adjacent groups. Therefore, a whole population can subdivide into subgroups with a given, average distance.
In the extreme case when $\gamma=0$, the whole population will not subdivide.
This automatic subdivision ability makes it
particularly suitable for highly nonlinear, multimodal optimisation problems.

In addition, the parameters in FA can be tuned to
control the randomness as iterations proceed, so that convergence can also be sped up by tuning these parameters.
These above advantages makes it flexible to deal with continuous problems, clustering and classifications,
and combinatorial optimisation as well.

As an example, let use use two functions to demonstrate the computational cost saved by FA. For details, please see the more extensive studies by Yang \cite{YangFA}.
For De Jong's function with $d=256$ dimensions
\be f(\x)=\sum_{i=1}^{256} x_i^2. \ee
Genetic algorithms required
$25412 \pm 1237$ evaluations to get an accuracy of $10^{-5}$ of the optimal solution,
while PSO needed $17040 \pm 1123$ evaluations. For FA, we achieved the same
accuracy by $5657 \pm 730$. This save about $78\%$ and $67\%$ computational cost, compared to GA and PSO, respectively.

For Yang's forest function
\be f(\x)=\Big(\sum_{i=1}^d |x_i| \Big) \exp \Big[-\sum_{i=1}^d \sin(x_i^2)\Big], \quad -2 \pi \le x_i \le 2\pi, \ee
GA required $37079 \pm 8920$ with a success rate of $88\%$ for $d=16$,
and PSO required $19725 \pm 3204$ with a success rate of $98\%$.
FA obtained a 100\% success rate with just $5152 \pm 2493$.
Compared with GA and PSO, FA saved about 86\% and 74\%, respectively, of overall computational efforts.

As an example for automatic subdivision, we now use the FA to find the
global maxima of the following function
\be f(\x)=\Big( \sum_{i=1}^d |x_i| \Big) \exp \Big( -\sum_{i=1}^d x_i^2 \Big),
\label{fun-equ-50} \ee
with the domain $-10 \le x_i \le 10$ for all $(i=1,2,...,d)$
where $d$ is the number of dimensions. This function has
multiple global optima \cite{YangFA2010}.  In the case of $d=2$, we have $4$ equal maxima
$f_* =1/\sqrt{e} \approx 0.6065$ at $(1/2,1/2)$, $(1/2,-1/2)$,
$(-1/2,1/2)$ and $(-1/2,-1/2)$ and a unique global minimum at $(0,0)$.

\begin{figure}
 \centerline{\includegraphics[width=4in,height=3in]{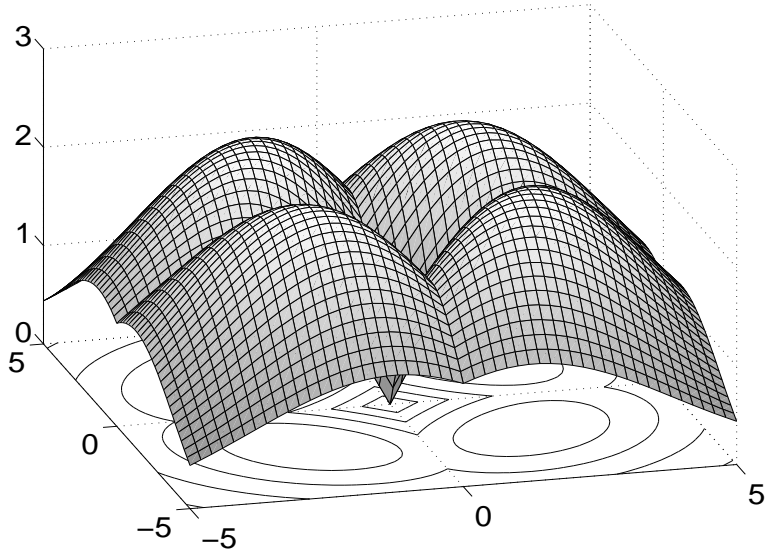} }
\caption{Four global maxima at $(\pm 1/2,\pm 1/2)$. \label{Yang-fig1} }
\end{figure}
\begin{figure}
 \centerline{\includegraphics[height=2.5in,width=3in]{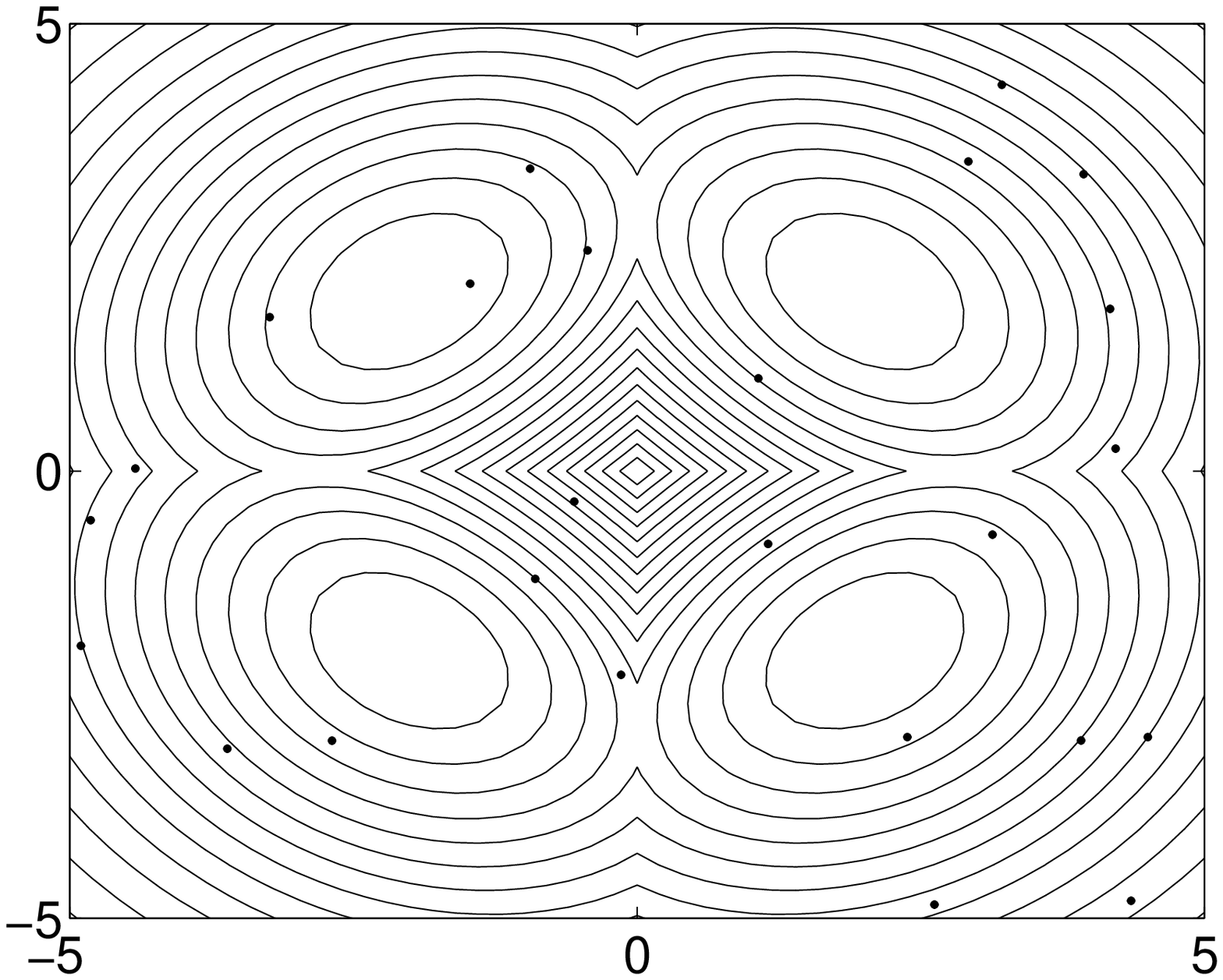}
 \includegraphics[height=2.5in,width=3in]{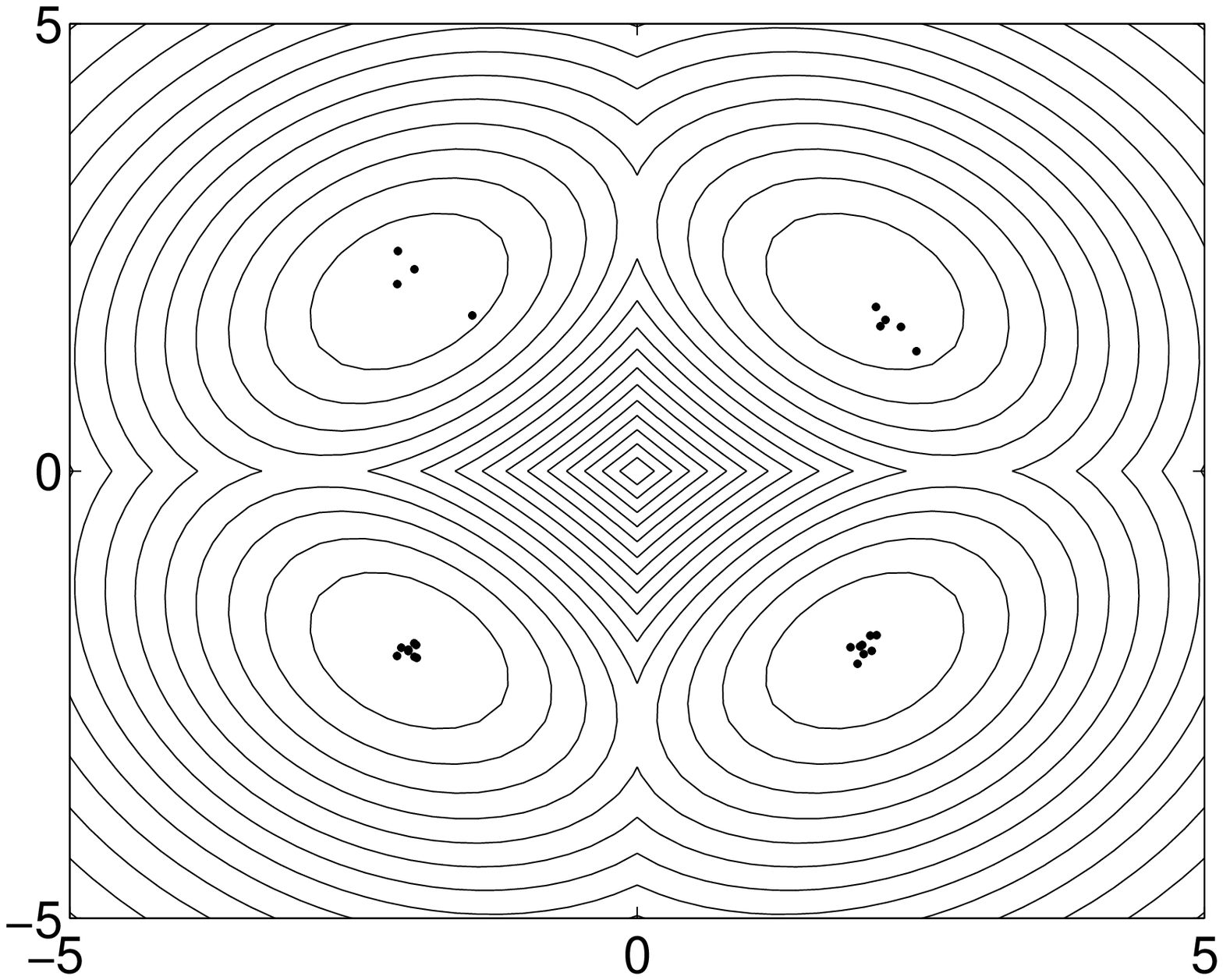} }
\vspace{-5pt}
\caption{Initial locations of 25 fireflies (left) and their
final locations after $20$ iterations (right). \label{Yang-fig2} }
\end{figure}

In the 2D case, we have a four-peak function are shown in Fig. \ref{Yang-fig1},
and these global maxima
can be found using the implemented Firefly Algorithms after about 500 function
evaluations. This corresponds to $25$ fireflies evolving for $20$ generations or
iterations. The initial locations of 25 fireflies
and their final locations after 20 iterations are shown in Fig. \ref{Yang-fig2}.
We can see that FA can indeed automatically subdivide its population into
subgroups.

\section{Search Optimality}

\subsection{Intensification versus Diversification}

The main components of any metaheuristic algorithms are: intensification
and diversification, or exploitation and exploration \cite{Blum,Yang12}.
Diversification means to generate diverse solutions so as to explore the search space
on the global scale, while intensification means to focus on the
search in a local region by exploiting the information that a current
good solution is found in this region. This is in combination
with the selection of the best solutions.

Exploration in metaheuristics can be achieved often by
the use of randomization \cite{Blum, Yang, YangFA},
which enables an algorithm to have the ability to jump out of any local optimum so as to explore the
search globally. Randomization can also be used for local search around the current best
if steps are limited to a local region. When the steps are large, randomization can
explore the search space on a global scale. Fine-tuning the right amount of randomness
and balancing local search and global search are crucially important
in controlling the performance of any metaheuristic algorithm.

Exploitation is the use of local knowledge of the search and solutions found so far
so that new search moves can concentrate on the local regions or neighborhood
where the optimality may be close; however, this local optimum may not be the
global optimality. Exploitation tends to use strong local information such as gradients,
the shape of the mode such as convexity, and the history of the search process.
A classic technique is the so-called hill-climbing which uses the local gradients
or derivatives intensively.

Empirical knowledge from observations and simulations of the convergence behaviour of
common optimisation algorithms suggests that exploitation tends to increase the
speed of convergence, while exploration tends to decrease the convergence rate of the
algorithm. On the other hand, too much exploration increases the probability of finding the global
optimality, while strong exploitation tends to make the algorithm being trapped in a
local optimum. Therefore, there is a fine balance between the right amount of exploration
and the right degree of exploitation. Despite its importance,
there is no practical guideline for this balance.

\subsection{Landscape-Based Optimality or Algorithm Optimality?}

It is worth pointing out that the balance for exploration and exploitation was often
discussed in the context of optimisation algorithms; however, this can be more often
related to the problem type of interest, and thus such balance is problem-specific, depending on the
actual landscape of the objective function. Consequently,  there may be no universality at all.
Therefore, we may have to distinguish landscape-dependent optimality for exploration and
exploitation, and algorithm-based optimality. The former focuses on landscape, while the latter focuses on algorithms.

Landscape-based optimality focuses on the information of the search landscape and it is hoped that
a better (or even best) approach/algorithm may find the optimal solutions with the minimum
effort (time, number of evaluations), while algorithm-based approach treats objective functions
as a black-box type and tries to use partly the available information during iterations to work out the
best ways to find optimal solutions. As to which approach is better, the answer may depend on one's
viewpoint and focus. In any case, a good combination of both approaches may be needed to
reach certain optimality.

\subsection{Intermittent Search Strategy}

Even there is no guideline in practice for landscape-based optimality, some preliminary
work on the very limited cases exists in the
literature and may provide some insight into the possible choice of parameters so as to balance these components.
In fact, intermittent search strategy is a landscape-based search strategy \cite{Ben}.

Intermittent search strategies concern an iterative strategy consisting of a slow phase and a fast phase \cite{Ben,Ben2}.
Here the slow phase is the detection phase by slowing down and intensive, static local search techniques,
while the fast phase is the search without detection and can be considered as an exploration technique.
For example, the static target detection with a small region of radius $a$ in a much larger region $b$
where $a \ll b$ can be modelled as a slow diffusive process in terms of random walks with a diffusion
coefficient $D$.

Let $\tau_a$ and $\tau_b$ be
the mean times spent in intensive detection stage and the time spent in the exploration stage,
respectively, in the 2D case. The diffusive search process is governed by
the mean first-passage time satisfying the following equations \cite{Ben2}
\be D \nabla_r^2 t_1 + \frac{1}{2 \pi \tau_a} \int_0^{2 \pi} [t_2(r) - t_1(r)] d\theta+1=0, \ee
\be u \cdot \nabla_r t_2(r) -\frac{1}{\tau_b} [r_2(r)-t_1(r)]+1=0, \ee
where $t_2$ and $t_1$ are mean first-passage times during the
search process, starting from slow and fast stages, respectively, and $u$ is the mean search speed.

After some lengthy mathematical analysis, the optimal balance of these two stages can be estimated as
\be r_{\rm optimal}=\frac{\tau_a}{\tau_b^2} \approx \frac{D}{a^2} \frac{1}{[2-\frac{1}{\ln(b/a)}]^2}. \label{equ-balance} \ee
Assuming that the search steps have a uniform velocity $u$ at each step on average, the minimum times required
for each phase can be estimated as
\be \tau_a^{\min} \approx \frac{D}{2 u^2} \frac{\ln^2(b/a)}{[2 \ln(b/a)-1]}, \ee
and
\be \tau_b^{\min} \approx \frac{a}{u} \sqrt{\ln(b/a)-\frac{1}{2}}. \ee
When $u \rightarrow \infty$, these relationships lead to the above optimal ratio of two stages.
It is worth pointing out that the above result is only valid for 2D cases, and there is no general results
for higher dimensions, except in some special 3D cases \cite{Ben}. Now let us use this limited results to help
choose the possible values of algorithm-dependent parameters in firefly algorithm \cite{Yang,YangFA}, as an example.

For higher-dimensional problems, no result exists. One possible extension is to
use extrapolation to get an estimate. Based on the results on 2D and 3D cases \cite{Ben2},
we can estimate that for any $d$-dimensional cases $d \ge 3$
\be \frac{\tau_1}{\tau_2^2} \sim O\Big(\frac{D}{a^2}\Big), \quad \tau_m \sim O\Big( \frac{b}{u} (\frac{b}{a})^{d-1} \Big), \ee
where $\tau_m$ the  mean search time or average number of iterations.
This extension may not be good news for higher dimensional problems, as the mean number of function evaluations to
find optimal solutions can increase exponentially as the dimensions increase. However, in practice,
we do not need to find the guaranteed global optimality, we may be satisfied with suboptimality,
and sometimes we may be `lucky' to find such global optimality even with a limited/fixed number of iterations.
This may indicate there is a huge gap between theoretical understanding and the observations as well as run-time
behaviour in practice. More studies are highly needed to address these important issues.

\section{Numerical Experiments}

\subsection{Landscape-Based Optimality: A 2D Example}

If we use the 2D simple, isotropic random walks for local exploration to demonstrate landscape-based optimality, then we have
\be D\approx \frac{s^2}{2}, \ee
where $s$ is the step length with a jump during a unit time interval or each iteration step.
From equation (\ref{equ-balance}), the optimal ratio of exploitation and exploration
in a special case of $b \approx 10 a$ becomes
\be \frac{\tau_a}{\tau_b^2} \approx 0.2. \ee
In case of $b/a \rightarrow \infty$, we have $\tau_a/\tau_b^2 \approx 1/8.$
which implies that more times should spend on the exploration stage. It is worth pointing out that
the naive guess of 50-50 probability in each stage is not the best choice. More efforts
should focus on the exploration so that the best solutions found by the algorithm can be
globally optimal with possibly the least computing effort. However, this case may be implicitly
linked to the implicit assumptions that the optimal solutions or search targets are
multimodal. Obviously, for a unimodal problem, once we know its modality, we should focus more on the
exploitation to get quick convergence.

In the case studies to be described below, we have used the firefly algorithm to find the
optimal solutions to the benchmarks. If set $\tau_b=1$ as the reference timescale, then
we found that the optimal ratio is between 0.15 to 0.24, which are roughly close to the above theoretical result.

\subsection{Standing-Wave Function}
Let us first use a multimodal test function to see how to find the fine balance between
exploration and exploitation in an algorithm for a given task. A standing-wave test function
can be a good example \cite{YangFA,YangBook}
\be f(\x)=1+ \Big\{\exp[-\sum_{i=1}^d (\frac{x_i}{\beta})^{10}]
 - 2 \exp[-\sum_{i=1}^d (x_i-\pi)^2]  \Big\} \cdot \prod_{i=1}^d \cos^2 x_i, \ee
which is multimodal with many local peaks and valleys. It has a unique global minimum at
$f_{\min}=0$ at $(\pi,\pi,...,\pi)$ in the domain $-20 \le x_i \le 20$ where $i=1,2,...,d$ and $\beta=15$.
In this case, we can estimate that $R=20$ and $a \approx \pi/2$, this means that $R/a \approx 12.7$, and
we have in the case of $d=2$
\be p_e \approx \tau_{\rm optimal} \approx \frac{1}{2 [2-1/\ln (R/a)]^2} \approx 0.19. \ee
This indicate that the algorithm should spend 80\% of its computational effort on global explorative search,
and 20\% of its effort on local intensive search.

For the firefly algorithm, we have used $n=15$ and $1000$ iterations.
We have calculated the fraction of iterations/function evaluations
for exploitation to exploration. That is, $Q=$ exploitation/exploration,
thus $Q$ may affect the quality of solutions. A set of 25 numerical experiments
have been carried out for each value of $Q$ and the results are summarized in Table 1.

\begin{table}
\begin{center}
\caption{Variations of $Q$ and its effect on the solution quality. }
\begin{tabular}{|l|l|l|l|l|l|l|l|l}
\hline
$Q \!\!$ & 0.4 & 0.3 & 0.2  & 0.1 & $0.05 \!\!$ \\ \hline
$f_{\min}\!\!$ & 9.4e-11 & 1.2e-12 & 2.9e-14 & 8.1e-12 & 9.2e-11 \\
\hline
\end{tabular}
\end{center}
\end{table}

This table clearly shows that $Q \approx 0.2$ provides the optimal balance of local exploitation
and global exploration, which is consistent with the theoretical estimation.

Though there is no direct analytical results for higher dimensions, we can expect that
more emphasis on global exploration is also true for higher dimensional optimisation problems.
Let us study this test function for various higher dimensions.

\subsection{Comparison for Higher Dimensions}

As the dimensions increase, we usually expect the number of iterations of finding the global
optimality should increase. In terms of mean search time/iterations,
B\'enichou et al.'s intermittent search theory suggests that \cite{Ben,Ben2}
\be \tau_m\Big|_{(d=1)} = \frac{2b}{u} \sqrt{\frac{b}{3a}}, \ee
\be \tau_m \Big|_{(d=2)} =\frac{2 b^2}{a u} \sqrt{\ln (b/a)}, \ee
\be \tau_m\Big|_{(d=3)}=\frac{2.2 b}{u} (\frac{b}{a})^2. \ee

For higher dimensions, we can only estimate the main trend based on the intermittent search strategy. That is,
\be \frac{\tau_1}{\tau_2^2} \sim O\Big(\frac{D}{a^2}\Big), \quad \tau_m \sim O\Big( \frac{b}{u} (\frac{b}{a})^{d-1} \Big), \ee
which means that number of iterations may increase exponentially with the dimension $d$.
It is worth pointing out that the optimal ratio between the two stage should be independent of the dimensions.
In other words, once  we find the optimal balance between exploration and exploitation, we can use the
algorithm for any high dimensions.

\begin{figure}
\centerline{\includegraphics[height=2.5in,width=3in]{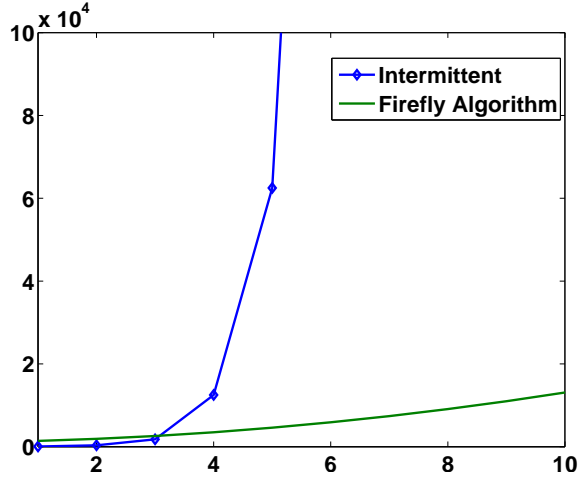}}
\caption{Comparison of the actual number of iterations with the theoretical results by the intermittent search strategy.
This clearly show that firefly algorithm is better than the intermittent search strategy. \label{fig-comp}}
\end{figure}

Now let us use firefly algorithm to carry out search in higher dimensions
for the above standing wave function and compare its performance with the implication of intermittent search strategy.
For the case of $b=20$, $a=\pi/2$ and $u=1$, Fig. \ref{fig-comp} shows the comparison of the numbers of iterations suggested by
intermittent search strategy and the actual numbers of iterations using firefly algorithm to obtain the
globally optimal solution with a tolerance or accuracy of 5 decimal places.
It can be seen clearly that the number of iterations needed by the intermittent search strategy increases
exponentially versus the number of dimensions, while the actual number of iterations used in the algorithm
only increases slightly, seemingly weakly a low-order polynomial. This suggests that firefly algorithm is very efficient
and requires far fewer (and often many orders lower) number of function evaluations.

\section{Conclusions}

Nature-inspired metaheuristic algorithms have gained popularity, which is partly due to
their ability of dealing with nonlinear global optimisation problems. We have reviewed the
fundamentals of firefly algorithm, the latest developments with diverse applications.
As the time of writing, a quick Google search suggests that there are about 323 papers on
firefly algorithms from 2008. This review can only cover a fraction of the literature.
There is no doubt that firefly algorithm will be applied in solving more challenging problems in the near future,
and its literature will continue to expand.

On the other hand, we have also highlighted the importance of exploitation and exploration
and their effect on the efficiency of an algorithm. Then, we use the
intermittent search strategy theory as a preliminary
basis for analyzing these key components and ways to find the possibly optimal settings
for algorithm-dependent parameters.

With such insight, we have used the firefly algorithm
to find this optimal balance, and confirmed that firefly algorithm can indeed
provide a good balance of exploitation and exploration.
We have also shown that
firefly algorithm requires far fewer function evaluations. However, the huge differences between
intermittent search theory and the behaviour of metaheuristics in practice also suggest
there is still a huge gap between our understanding of algorithms and the actual behaviour of
metaheuristics. More studies in metaheuristics are highly needed.

It is worth pointing out that there are two types of optimality here. One optimality concerns that
for a given algorithm what best types of problems it can solve. This is relatively easy to answer because
in principle we can test an algorithm by a wide range of problems and then select the best
type of the problems the algorithm of interest can solve. On other hand, the other optimality concerns
that for a given problem what best algorithm is to find the solutions efficiently. In principle,
we can compare a set of algorithms to solve the same optimisation problem and hope to find
the best algorithm(s). In reality, there may be no such algorithm at all, and all test algorithms
may not perform well. Search for new algorithms may take substantial research efforts.

The theoretical understanding of metaheuristics is still lacking behind.
In fact, there is a huge gap between theory and applications.
Though theory lags behind, applications in contrast are very diverse and active with
thousands of papers appearing each year. Furthermore, there is another huge gap
between small-scale problems and large-scale problems. As most published
studies have focused on small, toy problems, there is no guarantee that
the methodology that works well for such toy problems will work for
large-scale problems. All these issues still remain unresolved both
in theory and in practice.

As further research topics, most metaheuristic algorithms require good modifications
so as to solve combinatorial optimisation properly. Though with great interest and many extensive
studies, more studies are highly needed in the area of combinatorial optimisation using metaheuristic
algorithms. In addition, most current metaheuristic research has focused on
small scale problems, it will be extremely useful if further research can focus on
large-scale real-world applications.


\begin{thebibliography}{99.}

\bibitem{Absh}
 A. A. Abshouri, M. R. Meybodi and A. Bakhtiary, New firefly algorithm based on multiswarm and
 learning automata in dynamic environments, Third Int. Conference on Signal Processing Systems (ICSPS2011),
 Aug 27-28, Yantai, China, pp. 73-77 (2011).

\bibitem{Azad}
Sina K. Azad, Saeid K. Azad, Optimum Design of Structures Using an Improved Firefly Algorithm,
{\it International Journal of Optimisation in Civil Engineering}, {\bf 1}(2), 327-340(2011).


\bibitem{Apo} Apostolopoulos T. and Vlachos A., (2011). Application of the Firefly Algorithm for
Solving the Economic Emissions Load Dispatch Problem, International Journal of Combinatorics, Volume 2011,
    Article ID 523806. http://www.hindawi.com/journals/ijct/2011/523806.html



\bibitem{Banati}
H. Banati and M. Bajaj, Firefly based feature selection approach,
{\it Int. J. Computer Science Issues}, {\bf 8}(2), 473-480 (2011).

\bibitem{Bansal}
J. C. Bansal and K. Deep, Optimisation of directional overcurrent relay times
by particle swarm optimisation, in: Swarm Intelligence Symposium (SIS 2008),
IEEE Publication, pp. 1-7, (2008).

\bibitem{Basu}
B. Basu and G. K. Mahanti, Firefly and artificial bees colony algorithm for synthesis of scanned and broadside linear array antenna, {\it Progress in Electromagnetic Research B.},  {\bf 32}, 169-190 (2011).

\bibitem{Ben}
O. B\'enichou, C. Loverdo, M. Moreau and R. Voituriez, Two-dimensional intermittent search
processes: An alternative to L\'evy flight strategies, {\it Phys. Rev.}, E{\bf 74}, 020102(R), (2006).

\bibitem{Ben2}
O. B\'enichou, C. Loverdo, M. Moreau, and R. Voituriez,  Intermittent search strategies,
{\it Review of Modern Physics}, {\bf 83}, 81-129 (2011).


\bibitem{Blum}
C. Blum  and A. Roli,  Metaheuristics in combinatorial optimisation: Overview and
conceptural comparision, {\it ACM Comput. Surv.}, {Vol. 35}, 268-308 (2003).

\bibitem{Chatt}
A. Chatterjee, G. K. Mahanti, and A. Chatterjee,
Design of a fully digital controlled reconfigurable switched beam conconcentric ring array
antenna using firefly and particle swarm optimisation algorithm,
{\it Progress in Elelectromagnetic Research B.}, {\bf 36}, 113-131 (2012).

\bibitem{Coelho}
L. dos Santos Coelho, D. L. de Andrade Bernert, V. C. Mariani,
a chaotic firefly algorithm applied to reliability-redundancy optimisation,
in: {\it 2011 IEEE Congress on Evolutionary Computation (CEC'11)}, pp. 517-521 (2011).

\bibitem{Durkota}
 K. Durkota, Implementation of a discrete firefly algorithm for the QAP problem within the sage framework, BSc thesis, Czech Technical University, (2011).

\bibitem{Fara}
S. M. Farahani, A. A. Abshouri, B. Nasiri and M. R. Meybodi, A Gaussian firefly algorithm,
{\it Int. J. Machine Learning and Computing}, {\bf 1}(5), 448-453 (2011).

\bibitem{Fara2}
S. M. Farahani, B. Nasiri and M. R. Meybodi, A multiswarm based firefly algorithm
in dynamic environments, in:{\it Third Int. Conference on Signal Processing Systems (ICSPS2011)},
Aug 27-28, Yantai, China, pp. 68-72 (2011)

\bibitem{Fister}
I. Fister Jr, I. Fister, J. Brest, X. S. Yang, Memetic firefly algorithm for combinatorial optimisation,
in: {\it Bioinspired Optimisation Methods and Their Applications} (BIOMA2012) edited by
B. Filipi\v{c} and J. \v{S}ilc, 24-25 May 2012, Bohinj, Slovenia, pp. 75-86 (2012).


\bibitem{Floudas}
C. A. Floudas and P. M. Pardolos, {\it Encyclopedia of Optimisation},
2nd Edition, Springer (2009).


\bibitem{Gandomi}
A. H. Gandomi, X. S. Yang, and A. H. Alavi, Cuckoo search algorithm: a metaheuristic approach
to solve structural optimisation problems, {\it Engineering with Computers},  {\bf 27},
article DOI 10.1007/s00366-011-0241-y,  (2011).

\bibitem{Giann}
G. Giannakouris, V. Vassiliadis and G. Dounias, Experimental study on a
hybrid nature-inspired algorithm for financial portfolio optimisation, SETN 2010,
Lecture Notes in Artificial Intelligence (LNAI 6040), pp. 101-111 (2010).


\bibitem{Hass}
T. Hassanzadeh, H. Vojodi and A. M. E. Moghadam, An image segmentation
approach based on maximum variance intra-cluster method and firefly algorithm,
in: {\it Proc. of 7th Int. Conf. on Natural Computation (ICNC2011)}, pp. 1817-1821 (2011).


\bibitem{Horng}
M.-H. Horng, Y.-X. Lee, M.-C. Lee and R.-J. Liou, Firefly metaheuristic algorithm for training the
radial basis function network for data classification and disease diagnosis,
in: {\it Theory and New Applications of Swarm Intelligence} (Edited by R. Parpinelli and H. S. Lopes),
pp. 115-132 (2012).

\bibitem{Horng2}
M.-H. Horng, Vector quantization using the firefly algorithm for image compression,
{\it Expert Systems with Applications}, {\bf 39}, pp. 1078-1091 (2012).

\bibitem{Horng3}
M.-H. Horng and R.-J. Liou, Multilevel minimum cross entropy threshold selection based on
the firefly algorithm, {\it Expert Systems with Applications}, {\bf 38}, pp. 14805-14811 (2011).


\bibitem{Jati}
G. K. Jati and S. Suyanto, Evolutionary discrete firefly algorithm for travelling salesman problem,
ICAIS2011, Lecture Notes in Artificial Intelligence (LNAI 6943), pp.393-403 (2011).

\bibitem{Kennedy}
J. Kennedy and R. Eberhart, Particle swarm optimisation, in: {\it Proc.
of the IEEE Int. Conf. on Neural Networks}, Piscataway, NJ, pp. 1942-1948 (1995).

\bibitem{Nandy}
S. Nandy, P. P. Sarkar, A. Das, Analysis of nature-inspired firefly algorithm based back-propagation neural
network training, {\it Int. J. Computer Applications}, {\bf 43}(22), 8-16 (2012).

\bibitem{Palit}
S. Palit, S. Sinha, M. Molla, A. Khanra, M. Kule, A cryptanalytic attack on the knapsack cryptosystem
using binary Firefly algorithm, in: {\it 2nd Int. Conference on Computer and Communication Technology (ICCCT)},
15-17 Sept 2011, India, pp. 428-432 (2011).

\bibitem{Parp}
R. S. Parpinelli and H. S. Lopes, New inspirations in swarm intelligence: a survey,
{\it Int. J. Bio-Inspired Computation}, {\bf 3}, 1-16 (2011).


\bibitem{Rajini}
A. Rajini, V. K. David, A hybrid metaheuristic algorithm for classification using micro array data,
{\it Int. J. Scientific \& Engineering Research}, {\bf 3}(2), 1-9 (2012).

\bibitem{Ramp}
B. Rampriya, K. Mahadevan and S. Kannan, U
nit commitment in deregulated power system using Lagrangian firefly algorithm,
{\it Proc. of IEEE Int. Conf. on Communication Control and Computing Technologies (ICCCCT2010)},
pp. 389-393 (2010).

\bibitem{Sayadi} Sayadi M. K., Ramezanian R. and Ghaffari-Nasab N., (2010).
A discrete firefly meta-heuristic with local search for makespan minimization in permutation flow shop scheduling problems, Int. J. of Industrial Engineering Computations, {\bf 1}, 1--10.

\bibitem{Senthil}
J. Senthilnath, S. N. Omkar, V. Mani, Clustering using firely algorithm: performance study,
{\it Swarm and Evolutionary Computation}, {\bf 1}(3), 164-171 (2011).


\bibitem{Yang}
X. S. Yang, {\it Nature-Inspired Metaheuristic Algorithms}, Luniver Press, UK, (2008).

\bibitem{YangFA}
X. S. Yang, Firefly algorithms for multimodal optimisation,
{\it Proc. 5th Symposium on Stochastic Algorithms, Foundations and Applications},
(Eds. O. Watanabe and T. Zeugmann),  Lecture Notes in Computer Science, 5792: 169-178 (2009).

\bibitem{YangBook}
X. S. Yang, {\it Engineering Optimisation: An Introduction with Metaheuristic Applications},
John Wiley and Sons, USA (2010).

\bibitem{YangBA}
X. S. Yang,   A new metaheuristic bat-inspired algorithm,
in: {\it Nature Inspired Cooperative Strategies for Optimisation} (NICSO 2010)
(Eds. J. R. Gonzalez et al.), Springer, SCI {Vol. 284}, 65-74 (2010).

\bibitem{YangFA2010}
X. S. Yang,  Firefly algorithm, stochastic test functions and design optimisation,
{\it Int. J. Bio-Inspired Computation}, {\bf 2}(2), 78-84 (2010).

\bibitem{YangDeb}
X. S. Yang  and S. Deb,  Cuckoo search via L\'evy flights,
{\it Proceeings of World Congress on Nature \& Biologically Inspired
Computing} (NaBIC 2009, India), IEEE Publications, USA, pp. 210-214 (2009).

\bibitem{YangFAC}
X. S. Yang, Chaos-enhanced firefly algorithm with automatic parameter tuning,
{\it Int. J. Swarm Intelligence Research}, {\bf 2}(4), pp. 1-11 (2011).

\bibitem{Yang12}
X. S. Yang, Swarm-based metaheuristic algorithms and no-free-lunch theorems,
in: {\it Theory and New Applications of Swarm Intelligence} (Eds. R. Parpinelli and H. S. Lopes),
Intech Open Science, pp. 1-16 (2012).

\bibitem{YangAPSO}
X. S. Yang, S. Deb and S. Fong, (2011). Accelerated particle swarm optimization
and support vector machine for business optimization and applications,
Networked Digital Technologies (NDT'2011),
Communications in Computer and Information Science, Vol. 136, Part I,
pp. 53-66.

\bibitem{YangMOFA}
X. S. Yang, Multiobjective firefly algorithm for continuous optimization,
{\it Engineering with Computers}, Online Fist, DOI: 10.1007/s00366-012-0254-1 (2012).

\bibitem{Yousif}
A. Yousif, A. H. Abdullah, S. M. Nor, A. A. abdelaziz,
Scheduling jobs on grid computing using firefly algorithm,
{\it J. Theoretical and Applied Information Technology}, {\bf 33}(2), 155-164 (2011).

\bibitem{Zaman}
M. A. Zaman and M. A. Matin, Nonuniformly spaced linear antenna array design using firefly algorithm,
{\it Int. J. Microwave Science and Technology}, Vol. 2012, Article ID: 256759, (8 pages), 2012.
doi:10.1155/2012/256759

\end{thebibliography}
\end{document}